\documentclass[11]{article}
\usepackage{latexsym}
\usepackage{amssymb}
\usepackage{amsmath}
\usepackage[pctex32]{graphics}

\newcommand{\real}{\mathbb R}
\newcommand{\diam}{{\rm diam}}
\newtheorem{theorem}{Theorem}[section]
\newtheorem{corollary}{Corollary}[section]
\newtheorem{definition}{Definition}[section]
\newtheorem{lemma}{Lemma}[section]
\numberwithin{equation}{section}

\begin{document}
\title{\normalsize\bf NECESSARY AND SUFFICIENT CONDITION OF MORSE-SARD THEOREM FOR REAL VALUED FUNCTIONS
\footnotetext[0]{{\it Math Subject Classification} Primary 58K05; Secondary 58C25}
}
\author{\normalsize AZAT AINOULINE}
\date{February 23,2003}
\maketitle

\begin{abstract}
Necessary and sufficient condition is given for a set $A\subset {\mathbb R}^1$ to be a subset of the critical values set  for a $C^k$ function $f:{\mathbb R}^m\rightarrow{\mathbb R}^1$ .
\end{abstract}

\section*{Introduction}       
The well known Morse Theorem \cite{Morse} states that the critical values set $C_vf$ of a $C^k$ map $f:\real^m\rightarrow\real^1$ is of measure zero in $\real^1$~if $k\geqslant m$,~where $C_vf=f(C_pf)$~and $C_pf=\{x\in \real^m ; rank Df=0\}$~is the critical points set of $f$. Some generalizations of the necessary condition were made by a number of scientists. But the question is still open: whether a given measure zero set $A\subset R^1$ can be a set of critical values of some $C^k$ function $f:\real^m\rightarrow \real^1$?  In other words, what is necessary and sufficient condition of Morse theorem ?
A successful approach to describe the critical values sets was made by Yomdin \cite{Yomdin},\cite{Yomdin2},~Bates and Norton \cite{Norton-Bates} using a notion that in this paper we call $0_k$-sets:
\begin{definition}
\label{k degree set}
\rm  Let $A$ be a compact subset of $\real$. We define a countable set $Z(A)=\{(\alpha,\beta)\subseteq\real\setminus A~;~\alpha,\beta\in A,~\alpha<\beta \}$. We call $A$ a {\bf set of \boldmath $k$-degree, \unboldmath$k>0$}, if the series
$\sum_{z\in Z(A)} |z|^{\frac{1}{k}}~~~\text{converges},$~~and designate  $A$  as {\boldmath $0_k$\bf -set}~ if $A$ is set of $k$-degree~ with measure 0.
\end{definition}  
The following theorem was proven by Bates and Norton:

\noindent
{\bf Theorem (\cite{Norton-Bates}).}\it~~
A compact set $B=C_vf$~~for some function $f\in C^k(\real,\real)$~~with compactly supported derivative if and only if $B$~~is $0_k$-set.\rm\\

\noindent
And according to \cite{Norton-Bates}, results equivalent to the theorem above~~were obtained by Yomdin in his unpublished paper \cite{Yomdin2}, where he made a conjecture equivalent to the following:

\noindent
{\bf Conjecture.}\it~~ For $k\geqslant n$,~~a compact set $B=C_vf$~~for some function $f\in C^k(\real^n,\real)$~~with compactly supported derivative if and only if $B$~~is $0_{k/n}$-set.\\
\rm

In this paper author gives  necessary and sufficient condition for a set $A\subset \real^1$ to be a subset of the critical value set of a $C^{<k+\lambda}$ function $f:\real^m\rightarrow\real^1$.\\

\noindent
{\bf Main Theorem.}\it~~
A set $A\subseteq C_vf$ for some $C^{<k+\lambda}$ function \mbox{$f:\real^m\rightarrow \real^1$}, $m \leqslant k,$~~if and only if~~$A$ is subset of a $\sigma\!\!-\!\!0_{<\frac{k+\lambda}{m}}$ set.\\
\rm

\noindent
And as a corollary, the necessary and sufficient condition is given for a set $A\subset\real^n$ to be  subset of the image of the critical points set of rank zero of a $C^{<k+\lambda}$ function  $f:\real^m\rightarrow \real^n$.\\

\noindent
{\bf Corollary.}\it~~
A set $A\subseteq C_vf$~  for some $C^{<k+\lambda}$ function $f:\real^m\rightarrow\real^n,~~m\leqslant k,$~ if and only if its projection on any axis in $\real^n$ is a subset of a $\sigma\!-0_{<\frac{k+\lambda}{m}}$ set,\rm\\

\noindent
where we designate : \newline
$\bullet$~for $k\in [1,\infty)$~~a compact set $A\subseteq \real$~~as $0_{<k}$-set if $A$~~is $0_t$-set for every positive real $t<k$,\newline
$\bullet~A$ as a {\boldmath$\sigma\!-0_k$~\bf set}~({\boldmath$\sigma\!-0_{<k}$~\bf set})~~if $A=\bigcup _{i\in \mathbb N} A_i$, where $A_{i}$ is $0_k$-set ($0_{<k}$-set)~~$\forall~ i \in \mathbb N$.

\section{Definitions and Preliminary Lemmas}

For $k\in\mathbb N$, a map $f:\real^n\to\real^m$ is of class $C^k$ if it possesses continuous derivatives of all orders $\leqslant k$. For $\lambda\in[0,1)~f$ belongs to $C^{k+\lambda}$ provided $f\in C^k$ and the $k^{th}$ derivative $D^kf $~~satisfies a local H\"{o}lder condition: For each $x_0\in\real^n$, there exists a neighborhood $U$ of $x_0$ and a constant $M$ such that
\[\|D^kf(x)-D^kf(y)\|\leqslant M|x-y|^\lambda\]
for  all $x,y\in U$.~~If $f\in C^t$~~for all $t<k+\lambda$~~we write $f\in C^{<k+\lambda}.$ \\ 
\noindent
{\bf Generalized Morse Theorem.}\it\newline
Let $k,m\in\mathbb N$, and $\beta\in[0,1).$~ Let $A$~ be a subset of $\real^m$.\newline
a){\rm (see Norton~\cite{Norton})}~There exist subsets~ $G_i~~(i=0,1,2\dots)$ of $A $~~with~ $A=\bigcup_{i=1}^\infty G_i$ ~and any~ $f\in C^{k+\lambda}(\real^m,\real)$ ~~critical on~ $A$ ~satisfies for each~ $i~:~~|f(x)-f(y)|\leqslant M_i|x-y|^{k+\lambda}$
for all~ $x,y\in  G_i$ ~and some~~$M_i>0.$\newline
b)~There exist subsets~ $A_i~~(i=0,1,2\dots)$~ of $A$~~with~ $A=\bigcup_{i=1}^\infty A_i $~and any~ $f\in C^{<k+\lambda}(\real^m,\real)$~ critical on $A$ satisfies for each $i~:~~|f(x)-f(y)|\leqslant N_i|x-y|^{k+\lambda}$
for all~~$x,y\in  A_i$~and some~~$N_i>0.$\rm\\

\noindent
{\bf Proof.}\newline\noindent
a)See \cite{Norton}.\newline\noindent
b)~Similar to a).
\hspace*{\fill}$\Box$ 

\begin{definition}\rm For $m,n\in\mathbb N,~k\in\real$~a function $\psi:B\subseteq\real^m\rightarrow\real^n$~~is {\bf {\boldmath $D^k$}-function}~~ if $\exists K>0~~:~~~\forall b,b^\prime\in B~~~~~~|\psi(b)-\psi(b^\prime)|^k\leqslant K|b-b^\prime|.$
\end{definition}

\noindent
We can rewrite the {\bf Generalized Morse Theorem} in terms of $D^k$-functions:
\begin{corollary}\label{Generalized Morse-Sard Theorem}
Let $k,m\in\mathbb N$, and $\beta\in[0,1).$~ Let $A$~ be a subset of $\real^m$.\newline
a)~ There are subsets $G_i~~(i=0,1,2\dots)$~ of $~A~$ with~$ A=\bigcup_{i=1}^\infty G_i$~ and for any ~$f\in C^{k+\lambda}(\real^m,\real)$~ critical on~ $A,~~f\upharpoonright G_i$~ is a ~$D^{\frac{1}{k+\lambda}}$-function for each $i$.\newline
b)~ There are subsets $A_i~(i=0,1,2\dots)$~of~ $A$~ with $A=\bigcup_{i=1}^\infty A_i$~ and for any $f\in C^{k+\lambda}(\real^m,\real)$~critical on~$A,~~f\upharpoonright A_i $~is a  $D^{\beta}$-function for each $i$,~and every $\beta > \frac{1}{k+\lambda}.$
\end{corollary}

\noindent
{\bf Properties of $D^k$-functions:}\\
\noindent
{\bf 1) Extension on closure property.} If $f:A\subseteq\real^m\stackrel{D^k}{\rightarrow}\real^n$~ for some $k>0$,~ and $\overline A$~ is the closure of $A$,~ then  there exists and unique $C^0$ function $\overline{f}:\overline{A}\subseteq\real^m\to\real^n$~ such that $\overline{f}\upharpoonright A=f$. And $\overline{f}$~~is a $D^k$ function.\\

{\bf Proof.} $\forall a\in\overline{A}$ let us define $\overline{f}(a)=\lim\limits_{x\in A,x\to a}f(x)$, existence of the limit and the properties of $\overline{f}$ easy follow from the fact that $f\in D^k$.\\

\noindent
{\bf 2) Composition property.} If $g\in D^k$ and $f\in D^p$, then $g\circ f\in D^{kp}.$\\

\noindent
{\bf 3) Subsets property.} If $f:A\subseteq\real^m\stackrel{D^k}{\rightarrow}\real^m$~ for some $k>0$,~ then $f\upharpoonright B\in D^k$~ for any $B\subseteq A.$\\

Now we define a set ~$K^n_0 = \{Q_{i_0},~ i_0\in\mathbb N\}$,~~where every $Q_{i_0}$ is a closed cube in $\mathbb R^n$ with side length 1 and every coordinate of any vertex of $Q_{i_0}$ is an integer.
\newline 
In general, having constructed the cubes of $K^n_{s-1}$,  divide each $Q_{i_0,i_1,i_2,....,i_{s-1}}\in K^n_{s-1}$ into $2^n$ closed cubes of side $\frac{1}{2^s}$,  and let $K^{n}_{s}$ be the set of all those cubes. More precisely we will write \newline
$K^{n}_{s}=\{ Q_{i_0,i_1,i_2,...,i_{s-1},i_s}~~;~~ Q_{i_0,i_1,i_2,...,i_{s-1},i_s} \subseteq Q_{i_0,i_1,i_2,...i_{s-1}}\in K^n_{s-1}, 1\leqslant i_s\leqslant 2^n \}.$

\noindent 
We use here a family of  continuous space filling curves $f_n :~[0,1] \stackrel{onto}{\longrightarrow} [0,1]^n,~~ n\in\mathbb N$ with  special properties:
\begin{align}\label{cubes preserving 1}
&\text{if}~\alpha\subseteq [0,1]~\text{and}~\exists s\in\mathbb N:\alpha \in K^1_{n\cdot s}~\Longrightarrow~f_n(\alpha) \subseteq \delta~\text{for~some}~\delta\in K^n_s\\
&\text{if}~\delta\subseteq [0,1]^n~\text{and}~\exists s\in\mathbb N:\delta \in 
K^n_s~\Longrightarrow~f_n^{-1}(int(\delta)) \subseteq\alpha~\text{for~some}~\alpha\in K^1_{n\cdot s}
\label{cubes preserving 2}
\end{align}
where~~$int(\delta)$~~is~the~set~of~interior~points~of~~$\delta$.

\begin{definition}
\rm We call a function $f_n:~[0,1]\stackrel{onto}{\rightarrow}[0,1]^n$ with the properties \eqref{cubes preserving 1},\eqref{cubes preserving 2} {\bf cubes preserving}.
\end{definition}

\begin{theorem}[\cite{Ainouline1}]
\label{big lemma}
For every $n \in\mathbb N$ there exists a continuous cubes preserving function  $f_n:[0,1]\stackrel{onto}{\longrightarrow} [0,1]^n.$
\end{theorem}

\begin{lemma}
\label{SP is D^k}
Any continuous space-filling cubes preserving function $f_n:[0,1]\stackrel{\rm onto}{\rightarrow}[0,1]^n$~ is a $D^n$-function.
\end{lemma}

\noindent
{\bf Proof.} Let $a,b\in[0,1],~a<b$, then there exists $s_0\in \mathbb N$ such that
$\frac{1}{2^{n(s_0+1)}}\leqslant b-a\leqslant\frac{1}{2^{ns_0}}$.
 Then $[a,b]\subseteq\alpha^\prime\cup\alpha^{\prime\prime}$~ for some $\alpha^\prime,\alpha^{\prime\prime}\in K^1_{ns_0}$, such that $\alpha^\prime\cap\alpha^{\prime\prime}\not=\emptyset$.  From the definition of cubes preserving function it follows that $f_n(\alpha^\prime)\subseteq \delta^\prime,~f_n(\alpha^{\prime\prime}) \subseteq\delta^{\prime\prime}$~ for some $\delta^\prime,\delta^{\prime\prime}\in K^n_{s_0}$~with $\delta^\prime\cap\delta^{\prime\prime}\not = \emptyset$~ because $\alpha^\prime\cap\alpha^{\prime\prime}\not=\emptyset$.
Now using $\diam(\delta^\prime\cup\delta^{\prime\prime})\leqslant 2\sqrt{n}\cdot\frac{1}{2^{s_0}}$, we get:
$|f_n(b)-f_n(a)|\leqslant\diam ( f_n(\alpha^\prime\cup\alpha^{\prime\prime}))\leqslant\diam(\delta^\prime\cup\delta^{\prime\prime})\leqslant 2\sqrt{n}\frac{1}{2^{s_0}}$
or
\begin{align}
\label{a not = b}
&|f_n(b)-f_n(a)|^n\leqslant(2\sqrt{n}\frac{1}{2^{s_0}})^n=2^{2n}n^{\frac{n}{2}}\frac{1}{2^{n(s_0+1)}}\leqslant K|b-a|
\end{align}
where $K=2^{2n}\cdot n^\frac{n}{2}$.~In case $a=b$ we have
\begin{align}
\label{a = b}
&|f_n(b)-f_n(a)|^n=0\leqslant K\cdot 0=K\cdot|b-a|.
\end{align}
From (\ref{a not = b}) and (\ref{a = b}) it follows that $f_n$ is $D^n$-function on $[0,1].$
\hspace*{\fill}$\Box$

\begin{lemma}{\rm (see Corollary from Lemma 2 in \cite{Ainouline} and also combinational lemma in \cite{Norton-Bates})}\label{lemma 2 of [2]}
Let $f:[a,b]\rightarrow\real^1$~ be continuous with $A\subseteq[a,b]$~ compact and $B:=f(A)$. There exists an injective function $\gamma:Z(B)\rightarrow Z(A)$~ such that for each $z\in Z(B)$~ with (say) $\gamma(z)=(x,x^\prime)$,~ one has $z\subseteq(f(x),f(x^\prime))$.\rm
\end{lemma}

\begin{lemma}\label{f(A) is k degree set}
Let $f:[a,b]\subseteq\real^1\rightarrow\real^1$~ be a continuous function such that $f\upharpoonright A$~ is $D^k$- function for some closed $A\subseteq[a,b]$~ and $k>0$,~~then $f(A)$~ is $\frac{1}{k}$-degree set.
\end{lemma}

\noindent
{\bf Proof.}
By the Lemma \ref{lemma 2 of [2]}, there exists an injective function $\gamma:Z(f(A))\rightarrow Z(A)$~ such that for each $z\in Z(f(A))$ with (say) $\gamma(z)=(x,x^\prime)\in Z(A)$,~ one has $z\subseteq(f(x),f(x^\prime))$~then $|z|^k\leqslant|f(x)-f(x^\prime)|^k$~ for each $z\in Z(f(A))$ and $x,x^\prime\in A$~~$(\text{Note:~}(x,x^\prime)\in Z(A) \text{~means that~}x,x^\prime\in A).$  On the other hand for any $x,x^\prime\in A~~~~|f(x)-f(x^\prime)|^k\leqslant K|x-x^\prime|$~ for some $K>0$ because $f\upharpoonright A\in D^k$. Hence for each $z\in Z(f(A))$:
$|z|^k\leqslant |f(x)-f(x^\prime)|^k\leqslant K|x-x^\prime|=K|\gamma(z)|$,
where $\gamma(z)=(x,x^\prime)\in Z(A)$. And consequently 
$\sum_{z\in Z(f(A))} |z|^k\leqslant K\sum_{z\in Z(f(A))} |\gamma(z)|,$~~but 
$\sum_{z\in Z(f(A))} |\gamma(z)|\leqslant\sum_{z\in Z(A)} |z|$ 
~(because the function $\gamma$~is injective),~~then \newline
$\sum_{z\in Z(f(A))} |z|^k \leqslant K\sum_{z\in Z(A)}|z|\leqslant K|b-a|<\infty.$
So that the series $\sum_{z\in Z(f(A))} |z|^k$\newline converges and by the Definition \ref{k degree set}~~ $f(A)$~ is $\frac{1}{k}$-degree set. 
\hspace*{\fill}$\Box$ \\

\section{Proof of necessary condition of Main Theorem}
  
\begin{theorem}\label{proposition}.Let $n\leqslant k\in \mathbb N,~\lambda \in [0,1)$\newline
a)~ If $f\in C^{k+\lambda}(\real^n,\real^1)$,~then $C_v f~$~ is a $\sigma\!\!-\!\!0_\frac{k+\lambda}{n}$ set .\newline
b)~ If $F\in C^{<k+\lambda}(\real^n,\real^1)$,~then $C_v F~$~ is a $\sigma\!\!-\!\!0_{<\frac{k+\lambda}{n}}$ set .
\end{theorem}

\noindent
{\bf Proof.}\newline
a)~By the Corollary \ref{Generalized Morse-Sard Theorem}, there exist $G_i~(i=0,1,2\dots)$~ such   that for each $i~~~f\upharpoonright G_i$~ is $D^{\frac{1}{k+\lambda}}$-function. Then for each $i$~~$f\upharpoonright\overline G_i$~ is $D^{\frac{1}{k+\lambda}}$~~by the \textsl{"Extension on closure property of $D^k$-functions"} and $C_pf=\cup^\infty_{i=1} \overline G_i$~~becouse $C_pf$~~is closed .

We may suppose without loss of generality that for each  $i~~~G_i$~ and, therefore, $\overline G_i$, is contained in some closed cube $C_i\subseteq\real^n$~ of side length 1 (because otherwise $G_i$ can be represented as countable union of sets contained in such cubes). 

Let for each $i\in \mathbb N$~~$g_i:[0,1]\stackrel{onto}{\rightarrow} C_i$~~be a continuous $D^n$ function that exists  by the Theorem \ref{big lemma} and Lemma \ref{SP is D^k}.Then on the closed set $g_i^{-1}(\overline G_i)~~g_i\upharpoonright(g_i^{-1}(\overline G_i)) \text{~is~}D^n-\text{function ~} $~ by the \textsl{"Subsets property of $D^k$-functions"}.
And by the \textsl{"Composition property of $D^k$-functions"}~~ $f\circ g_i\upharpoonright(g_i^{-1}(\overline G_i))$~is a $D^{\frac{n}{k+\lambda}}$-function .~So that by the Lemma \ref{f(A) is k degree set} and by the fact that $f\circ g_i$~ is continuous on $[0,1]$~ as a composition of continuous functions, it follows that the closed set
$f(g_i(g_i^{-1}(\overline G_i))) \text{~~is of~~}\frac{k+\lambda}{n}-\text{degree set},$
and therefore $f(\overline G_i)$~ is of $\frac{k+\lambda}{n}$-degree set.

By the Morse Theorem \cite{Morse} $f(C_pf)=C_vf$~ is a measure zero set.And becouse for each $i~~\overline G_i \subseteq C_pf$~~ it follows that $f(\overline G_i)$~ is $0_{\frac{k+\lambda}{n}}$-set ,  and consequently $f(C_pf)=C_vf$~~is a $\sigma\!-0_\frac{k+\lambda}{n}$~~set (recall that $C_pf=\bigcup_{i=0}^\infty \overline G_i$).

b)~Similar to a).
    The Theorem \ref{proposition} is proven, and thereby the prove of the necessary condition of the Main Theorem has been finished.

\section{Proof of sufficient condition of Main Theorem}

\begin{theorem}\label{theorem 3.1}
If $A$ is a $\sigma\!-0_{<s}$-set~ $(1<s\in \real)$,~then for every $n\in \mathbb N~~A\subseteq C_vF$ for some $C^{<sn}$ function $F:\real^n\rightarrow \real^1$.
\end{theorem}
\noindent
{\bf Proof.}
First let us fix such $s,n$~and prove that for any $0_{<s}$~set $B$~there exists a $C^{<sn}$~function $f:[0,1]^n\rightarrow \real$~with $B\subseteq C_vF.$  
Since $B$ is an $0_{<s}$-set,~ then the series $\sum_{z_n\in Z(B)}|z_n|^{1/t}$~~converges for any positive $t<s$. The case $B=\emptyset$ is trivial, so we suppose $B\not =\emptyset$.

  Let $P\geqslant 1$ denote the greatest integer less than $sn$~and $\{s_m,~m\in \mathbb N\}$~~ be a non-decreasing sequence such that $(P+sn)/2<s_m<sn,~~\lim\limits_{m\to\infty}s_m = sn$ ~~ and~~$\sum_{z_m\in Z(B)} |z_m|^{1/s_m}$~~ is convergent.  The existence of such sequence~$\{s_m\}$ can be seen from the fact that~~$\forall i\in\mathbb N,~~ \forall t\in \real^+,~~t<sn,$~~ there exists $m_i\in \mathbb N$~~such that $\sum_{m>m_i,z_m\in Z(B)}|z_m|^{1/t}<\frac{1}{2^i}$,~~ and $\sum^\infty_{i=1}\frac{1}{2^i}= 2$.

  For any closed $B^*\subseteq B$~let us designate the sum $\sum_{z_m\in Z(B^*)} |z_m|^{1/s_m} = G(B^*)$.  Further we construct the function $f$ following the scheme of Bates (\cite{Bates},sec.1.4) and replacing   definitions of $P,~R(a_k)$ and $\sigma_t(k)$ in his construction.

\subsection{Construction of $f$}

For $\beta\in(0,1/2)$, we can first define the following method for constructing $2^n$ cubes within any cube in $\real^n:$\newline
Supposing $Q\subset\real^n$ is a cube of side $L>0$ defined by the inequalities $|x_i-c_i|\leqslant L/2$, we specify $2^n$ subcubes within $Q$ with the inequalities $|x_i-c_i\pm L/4|\leqslant\beta L/2.$ Note that each subcube is separated from all other subcubes and the boundary of $Q$ by a distance $\geqslant L(1/4-\beta/2).$

\subsubsection{A Cantor set in $\real^n.$}

Let $Q_0\subset\real^n$ be the cube defined by $|x_i|\leqslant 1/2.$ For $i\in\mathbb Z^+$ let $\beta_i=(1/2)\cdot e^{-1/i}$ and set $\pi_k=(16k)^{-1}\cdot\prod_{i=1}^k\beta_i$. For $k\in\mathbb Z^+$ the index $a_k$ always denotes a $k$-tuple of numbers in $\{1,2,3,...,2^n\}$. 
We now define a system of subcubes in $Q_0$ as follows:\newline
\indent
a. Let $\{Q(a_1)\}$ be the $2^n$ subcubes constructed in $Q_0$ by the method described above with $\alpha=\beta_1$.\newline
\indent
b. Having defined $\{Q(a_k)\}$, construct $2^n$ subcubes within each $Q(a_k)$ according to the above process with $\alpha=\beta_{k+1}$ and label these subcubes $\{Q(a_k,i)\}$ for $i=1,...,2^n.$

Evidently, the cube $Q(a_k)$ has side length $L_k=\prod_{i=1}^k\beta_i;$ for integers $k\leqslant l$ the boundaries of distinct cubes $Q(a_k)$ and $Q(a_l^\prime)$ are separated by a distance $\geqslant(1/4-\beta_{k+1}/2)L_k\geqslant \pi_k.$

Let $\zeta$ denote the Cantor set defined by the cubes, i.e. the set of points in $Q_0$ contained in indefinitely many of subcubes constructed above.

\subsubsection{Mapping of $\zeta$.}

Let $R_0=B$. For each $k\in \mathbb Z^+$, choose decomposition of $R_0$~into $2^{nk}$ non-empty  closed set $\{R(a_k)\}$~such that $G(R(a_k))\leqslant M\cdot 2^{-nk}$~and $R(a_k)=\cup^{2^n}_{i=1}R(a_k,i)$. For each $a_k$, fix a point $r(a_k)\in R(a_k)$.

Define the map $f$~on $\zeta$~by the requirement that for each index $a_k, f(\zeta \cap Q(a_k))\subseteq R(a_k).$~The conditions imposed on the sets $R(a_k)$~then insure that $R_0\subseteq f(\zeta)$.

\subsubsection{Extension of $f$.}

Consider the cube $Q=Q(a_k)$ and its subcubes $Q_i=Q(a_k,i)$. For each $i=1,...,2^n$ choose a function $h_i:\real^n\rightarrow\real$  such that\newline\indent
(1) $h_i=1$ on a neighborhood of $Q_i;$\newline\indent
(2) $Supp(h_i)\subset Int~Q$, and $Supp(h_i)\cap Supp(h_j)=\emptyset$ whenever $i\not=j$.\newline
In view of condition (2) and the distance between the $Q_i$, we can choose the $h_i$ so that, for each $p\in\mathbb Z^+$,\newline\indent
(3) $\|D^ph_i\|\leqslant M_p(\pi_{k+1})^{-p}.$

Now we define the partial extension of $f$ to the region $Q\setminus\bigcup Q_i$ by 
\[f=r+\sum\limits_{i=1}^{2^n}(r_i-r)h_i,\]
where $r=r(a_k),~r_i=r(a_k,i).$

\subsubsection{Smoothness of $f$.}

 For $t\in\real$ define $\sigma_t(k)=2^{-(sn+t)k/2}\cdot(\pi_{k+1})^{-t}$. Since $\beta_k\to 1/2$ as $k\to\infty$, a simple estimation shows that $\pi_k$ is bounded below by $M^\prime\cdot2^{-k}k^{-2}$. Consequently, for $t<sn \quad \sigma_t(k)\to 0$ as $k\to\infty.$

To determine the smoothness of $f$, we first observe that $(P+sn)/2n<s$,~also if $k\geqslant 1$ is the largest integer such that $x,x^\prime\in Q(a_k)\cap\zeta,$ then $|x-x^\prime|\geqslant \pi_{k+1}$, and by the definition of $f$,
\begin{align*}
&|f(x)-f(x^\prime)|\leqslant diam(R(a_k))=\sum\nolimits_{z_m\in Z(R(a_k))} |z_m|\\
&\leqslant \left(\sum\nolimits_{z_m\in Z(R(a_k))} |z_m|^{2n/(P+sn)}\right)^{(P+sn)/2n} \leqslant (G(R(a_k)))^{(P+sn)/2n}\\
&\leqslant M^\prime \sigma_P(k)\pi_{k+1}\leqslant M^\prime\sigma_P(k)|x-x^\prime|^P.
\end{align*}
This implies that $f\upharpoonright\zeta$ is continuous; by condition (1), it follows that the partial extension defined above comprise a continuous extension of $f\upharpoonright\zeta$. Since condition (2) implies $\|D^pf\|=0$ on boundaries $\partial Q(a_k)$ for all $k,p\in\mathbb Z^+$, it follows furthermore that $f$ is $C^\infty$ on $\real^n\setminus\zeta$.

Condition (3) above implies that on $Q(a_k)\setminus\bigcup Q(a_k,i)$
\[\|D^pf\|\leqslant M_p(\pi_{k+1})^{-p}\cdot diam(R(a_k))\leqslant M^\prime_p\sigma_p(k).\]
Consequently $D^pf\to 0$ on approach to $\zeta$ whenever $p\leqslant P$. By an application of the mean-value theorem to the inequality found for $f\upharpoonright\zeta$ above, it follows that $D^pf=0$ on $\zeta$ for $p\leqslant P$, and so $f$ is $C^P$ on $\real^n.$

Given $t\in(P,sn)$, we observe that if $x\in\zeta$ and $k\gg 1$ is again the largest integer such that $x,y\in Q(a_k)$, then the same argument used above shows that
\[\|D^Pf(x)-D^Pf(y)\|/|x-y|^{t-P}=\|D^Pf(y)\|/|x-y|^{t-P}\leqslant M^{\prime\prime}_P\sigma_t(k+1).\]
Since $f\in C^\infty$ outside $\zeta$, this inequality implies $f\in C^t$ throughout $\real^n.$
Evidently, $rank Df=0$ on $\zeta$, and $B\subseteq C_vf$.\\

By construction, $f$ is constant outside $Q_0$. Now for a $\sigma\!-O_{<s}$ set $A$, we can define $C^{<sn}$ function  $F:\real^n\to\real^1$, such that $A\subseteq C_vF$:\newline
We take $\{Q^i_0\}_{i\in\mathbb N}$ - a discrete family of cubes of side length 1 in $\real^n$ and a family of $C^{<sn}$ functions $\{f_i:Q^i_0\to\real^1~;~A_i\subseteq C_vf_i\}$. For each $i$ we define $F=f_i\upharpoonright Q^i_0$. And we can choose distances between cubes allowing $F$ be $C^\infty$ on $\real^n\setminus\cup Q^i_0$, and $C^{<sn}$ on $\real^n.$
\hspace*{\fill}$\Box$ \\


\begin{thebibliography}{10}
\bibliographystyle{plain}
 \bibitem{Ainouline1}Ainouline, A.,Vector valued functions not constant on connected sets of critical points, 2003, ArXiv: math. GT/0404404.
 \bibitem{Ainouline}Ainouline, A. and Zvengrowski, P.,Critical Values of differential functions on the Reals, University of Calgary, Research Paper No.813,July 2001.
 \bibitem{Bates}Bates, S.M., On the image size of singular maps II, Duke Math.J. Vol.68, No.3, December, 1992, 463-476. 
 \bibitem{Morse}Morse, A.P.,The behavior of a function on its critical set,Annals of Math. 40(1939), 62-70.
 \bibitem{Norton}Norton, A., A critical set with nonull image has large Hausdorff dimension, Trans.Amer.Math.Soc., 296, No.1, July 1986, 367-376.
 \bibitem{Norton-Bates}Bates, S.M. and Norton, A., On sets of critical values in the real line,Duke Math.J. 83(1996), No.2, 399-413.
 \bibitem{Yomdin}Yomdin, Y.,The geometry of critical and near-critical values of differentiable mappings, Math.Ann. 264(1983), 495-515.
 \bibitem{Yomdin2}Yomdin, Y.,$\beta$-spread of sets in metric spaces and critical values of smooth functions, preprint, Max-Planck-inst.,Bonn, 1982.
\end{thebibliography}
\end{document}